# FORMS OF REPRESENTATION OF INTERPOLATION TRIGONOMETRIC SPLINES

**Denysiuk V.P.**, Doctor of Physical and Mathematical Sciences, Professor,

**Rybachuk L.V.**, Candidate (Ph. D.) of Physical and Mathematical Sciences, Docent

National Technical University of Ukraine "Igor Sikorsky Kyiv Polytechnic Institute", Ukraine
(kvomden@nau.edu.ua, rybachuk.liudmyla@lll.kpi.ua)

**Abstract**

Three forms of representation of trigonometric interpolation splines are considered, in particular, the representation by the coefficients of the interpolation trigonometric polynomial, the representation by trigonometric B-splines, which are considered in more detail, and the representation by fundamental trigonometric splines. The first and third forms of representation are generalized to the case of non-periodic functions.

## Introduction

The approximation, or representation, of an arbitrary known or unknown function through the set of some special functions can be considered as the central theme of the analysis. The term "special functions" refers to the classes of algebraic and trigonometric polynomials and their modifications; we assume that the classes of trigonometric polynomials include trigonometric series. As a rule, such special functions are easy to compute and have interesting analytical properties [1].

One of the most successful modifications of algebraic polynomials is polynomial splines, which are sewn from segments of these polynomials according to a certain scheme. The theory of polynomial splines has appeared relatively recently and is well developed (see, for example, [2], [3], [4] [5], etc.). The advantages of polynomial splines include their approximation properties [6]. The main disadvantage of polynomial splines is their piecewise structure, which greatly complicates their use in analytical transformations.

An important role in the theory of polynomial splines is played by normalized B-splines (hereinafter simply B-splines), which are a basis in the space of polynomial splines [3] [4]. The representation of polynomial splines through B-splines is used in many problems requiring analytical transformations, in particular, in solving differential and integral equations (see, for example, [7], [8]), etc.

Subsequently, it turned out [9], [10] that there are also modifications of trigonometric series that depend on several parameters and have the same properties as polynomial splines [11], [12], [12]; moreover, the class of such modified series is quite broad and includes the class of polynomial periodic splines. This gave rise to the name of the class of such series as trigonometric splines.

The important advantages of trigonometric splines are that they are represented by a single analytical expression and are carried out in the same way, regardless of their order; their construction does not depend on the order of the spline; the construction of polynomial splines of high degrees is associated with significant difficulties. The advantages of trigonometric splines include the fact that all the results of the theory of approximations obtained for periodic polynomial splines are transferred to these splines.

Unlike polynomial splines, trigonometric splines are represented by analytical expressions and have several forms of representation.

Therefore, in our opinion, it is relevant to study the representation of trigonometric splines through the previously identified trigonometric B-splines, as well as to consider other forms of representation of trigonometric splines.

**Purpose of the work**

The study of the representation of trigonometric interpolation splines through trigonometric B-splines and consideration of other forms of representation of trigonometric interpolation splines.

**Main text**

**1.** Let $C^k_{[0,\pi]}$ be a class of functions that have a continuous derivative of order $k$ on $[0,\pi]$. We will also consider the class $C^k_{[0,2\pi)}$ of $2\pi$-periodic functions with a continuous derivative of order $k$. Note that the classes $C^0_{[0,\pi]}$ and $C^0_{[0,2\pi)}$ contain continuous functions, and the classes $C^{-1}_{[0,\pi]}$ and $C^{-1}_{[0,2\pi)}$ contain piecewise continuous functions with discontinuities of the first kind.

Let on a segment $[0,2\pi)$ we have the grids $\Delta 1_N^{(I)} = \{x_j^{(I)}\}_{j=1}^N$, where $I$ is the grid indicator, $(I=0,1)$, $x_j^{(0)} = \frac{2\pi}{N}(j-1)$, $x_j^{(1)} = \frac{\pi}{N}(2j-1)$, $N = 2n+1$, ($n=1,2,...$). Similarly, on $[0,\pi]$ we define the uniform grids $\Delta 2_N^{(I)} = \{x_j^{(I)}\}_{j=1}^N$, where $x_j^{(0)} = \frac{\pi}{N-1}(j-1)$, $x_j^{(1)} = \frac{\pi}{2N}(2j-1)$, and $\Delta 3_N^{(I)} = \{x_j^{(I)}\}_{j=1}^N$, where $x_j^{(0)} = \frac{\pi}{N+1}j$, $x_j^{(1)} = \frac{\pi}{2N}(2j-1)$, $I$ is the grid indicator, ($I=0,1$), $N = 1,2,3,...$. We will not introduce separate notations for the nodes of these grids, since in each case we will determine which grid we are talking about. Note that all of the introduced meshes $\Delta 1_N^{(I)}$, $\Delta 2_N^{(I)}$, $\Delta 3_N^{(I)}$ are uniform; it is the uniform meshes that we will consider in this paper.

In [9]-[13], the Riemann convergence multiplier was studied (see, e.g., [14] [15] [16]), which has the form

$$\sigma(r,k) = \left(\frac{\sin\frac{\pi}{N}k}{\frac{\pi}{N}k}\right)^{1+r}. \quad (1)$$

In this paper, we will focus on the consideration of the familiar multiplier

$$\sigma 1(r,k) = \left(\frac{1}{k}\right)^{1+r}. \quad (2)$$

Note that the familiar constant convergence factors of (2) allow for broad generalizations, since arbitrary familiar infinitesimal factors equivalent to infinitesimal $\left(\frac{1}{k}\right)^{1+r}$ can be used as such factors. It is clear that other types of trigonometric splines can be obtained.

**2. Representation of trigonometric splines through the coefficients of interpolation trigonometric polynomials.** In this presentation, we will consider the grid $\Delta 1_N^{(I)}$ given on $[0,2\pi)$.

As is well known (see, e.g., [17]), a trigonometric polynomial $T_n(x)$ of order $n$ ($n=1,2,...$) of the form

$$T_n^{(I)}(x) = \frac{a_0^{(I)}}{2} + \sum_{k=1}^n \left[a_k^{(I)} \cos kx + b_k^{(I)} \sin kx\right], \quad (3)$$

interpolating the function $f(x)$ in the grid nodes $\Delta 1_N^{(I)}$, has the coefficients

$$a_0^{(I)} = \frac{2}{N}\sum_{j=1}^N f_j^{(I)}, \quad a_k^{(I)} = \frac{2}{N}\sum_{j=1}^N f_j^{(I)} \cos kx_j^{(I)}, \quad b_k^{(I)} = \frac{2}{N}\sum_{j=1}^N f_j^{(I)} \sin kx_j^{(I)}, \quad k=1,2,...,n. \quad (4)$$

The trigonometric $2\pi$-periodic interpolation spline $St(I1,I2,\sigma 1,r,q,t)$ with a crosslinking grid $I1$, interpolation grid $I2$, of order $r$, with a sign-constant convergence factor $\sigma 1$, and its derivatives of order $q$ ($q=1,2,...,r$), can be obtained by the following formulas

$$St(I1,I2,\sigma 1,r,q,t) = \frac{a_0}{2} I(q) + \tag{5}$$

$$+ \sum_{k=1}^{n} \left[ \left( \frac{C(I1,r,q,k,t)}{H(I1,I2,\sigma 1,r,k)} \cdot a_k^{1-I2} a1_k^{I2} + \frac{S(I1,r,q,k,t)}{H(I1,I2,\sigma 1,r,k)} \cdot b_k^{1-I2} b1_k^{I2} \right) \right],$$

where

$$I(q) = \begin{cases} 1, & \text{if } q=0, \\ 0 & \text{otherwise} \end{cases}; \tag{6}$$

$$H(I1,I2,\sigma 1,r,k) = \sigma 1(r,k) + \sum_{m+1}^{\infty} \left[ (-1)^{m(I1+I2)} \left( \sigma 1(r,mN+k) + (-1)^{1+r} \sigma 1(r,mN-k) \right) \right]; \tag{7}$$

$$C(I1,r,q,r,t) = \sigma 1(r,k) k^q \cos(kt + q\frac{\pi}{2}) +$$

$$+ \sum_{m=1}^{\infty} \left[ (-1)^{m \cdot I1} \left[ \sigma 1(r,mN+k)(mN+k)^q \cos\left( (mN+k)t + q\frac{\pi}{2} \right) + \right. \right. \tag{8}$$

$$\left. \left. + (-1)^{1+r} \sigma 1(r,mN-k)(mN-k)^q \cos\left( (mN-k)t + q\frac{\pi}{2} \right) \right] \right];$$

$$S(I1,r,q,r,t) = \sigma 1(r,k) k^q \sin(kt + q\frac{\pi}{2}) +$$

$$+ \sum_{m=1}^{\infty} \left[ (-1)^{m \cdot I1} \left[ \sigma 1(r,mN+k)(mN+k)^q \sin\left( (mN+k)t + q\frac{\pi}{2} \right) - \right. \right. \tag{9}$$

$$\left. \left. - (-1)^{1+r} \sigma 1(r,mN-k)(mN-k)^q \sin\left( (mN-k)t + q\frac{\pi}{2} \right) \right] \right].$$

It is clear that the spline $St(I1,I2,\sigma 1,r,0,t)$ is actually a Fourier series whose terms have the descending order $O(k^{-(1+r)})$; it follows that this spline is $2\pi$-periodic function and $St(I1,I2,\sigma 1,r,0,t) \in C_{[0,2\pi)}^{r-1}$. The derivatives of this spline of order $q$ $(q=1,2,...,r)$ are defined by formula (5) and $St(I1,I2,\sigma 1,r,q,t) \in C_{[0,2\pi)}^{r-1-q}$. Note that the splines $St(0,0,\sigma,r,0,t)$ at odd values of the parameter $r$ coincide with the polynomial interpolation simple periodic spline on $[0,2\pi)$ and its derivatives of order $q$ $(q=1,2,...,r)$.

Note that the functions $\frac{C(I1,r,q,k,t)}{H(I1,I2,\sigma 1,r,k)}$ and $\frac{S(I1,r,q,k,t)}{H(I1,I2,\sigma 1,r,k)}$, calculated by formulas (7)-(9), form a system of functions orthogonal to the set of discrete equidistant points and belong to the class $C_{[0,2\pi)}^{r-1}$. With this approach, the coefficients $a_0^{(I)}$, $a_k^{(I)}$, $b_k^{(I)}$ ($k=1,2,...,n$) can be considered as coefficients of discrete (finite) Fourier series [18], [19].

Trigonometric splines $St(I1,I2,\sigma 1,r,q,t)$ can be generalized in several ways. For example, smoothing methods [19] or similar methods $\lambda$-summation (see, e.g., [16]), which are well known in the theory of trigonometric series, can be applied to such splines. These methods differ in that smoothing is applied to sums with approximate Fourier coefficients (or to finite Fourier series), and $\lambda$-summation methods are applied to partial sums with exact Fourier coefficients (or to partial sums of classical Fourier series); in our opinion, in many cases in practice this difference can be neglected and the factors of the theory of $\lambda$-summation can be used in smoothing problems. In this approach, depending on the chosen $\lambda$-summation factors, trigonometric splines of Fejér, Vallée-Poussin, Hamming-Tukey, Rogosinski, etc. are obtained. Another direction of generalization is the method of introducing factors that change the ratio between the low-frequency, medium-frequency, and high-frequency components of trigonometric splines [11].

If the function $f(x)$ is not periodic, it is advisable to use incomplete (even or odd) splines. We will begin our consideration of such splines by considering even splines, which will be considered on the grids $\Delta 2_N^{(I)}$, given by $[0, \pi]$.

First of all, let's note that the coefficients of an even trigonometric polynomial interpolating the function $f(x)$ on the grid $\Delta 2_N^{(0)}$ are calculated using the following formulas:

$$a0 = \frac{2}{N-1}\left[.5(f_1 + f_N) + \sum_{j=2}^{N-1} f_j\right] \tag{10}$$

$$a_k = \frac{2}{N-1}\left[.5(f_1 + f_N \cdot \cos(kx_N)) + \sum_{j=2}^{N-1} f_j \cdot \cos(kx_j)\right], \tag{10a}$$

The even spline $Stc(0,0,\sigma 1,r,q,t)$ and its derivatives of the order $q$ $(q=1,2,...,r)$, are calculated by the formulas

$$Stc(0,0,\sigma 1,r,q,t) = .5a_0 I(q) + \sum_{k=1}^{N-2} a_k \frac{C(\sigma 1,r,q,k,t)}{H(\sigma 1,r,k)} + .5a_{N-1}\frac{C(\sigma 1,r,q,N-1,t)}{H(\sigma 1,r,N-1,k)}, \tag{11}$$

where $I(q)$ is determined by formula (5);

$$H(\sigma 1, r, k) = \sigma 1(r, k) + \sum_{m=1}^{\infty}\left[\sigma 1(r, 2m(N-1)+k) + \sigma 1(r, 2m(N-1)-k)\right], \tag{12}$$

$$C(\sigma 1, r, q, k, t) = \sigma 1(r,k)k^q \cos(kt + q\frac{\pi}{2}) +$$
$$+ \sum_{m=1}^{\infty}\left[\sigma 1((r, 2m(N-1)+k))((2m(N-1)+k)^q \cos((2m(N-1)+k)t + q\frac{\pi}{2}) + \right. \tag{13}$$
$$\left. + \sigma 1((r, 2m(N-1)-k))(2m(N-1)-k)^q \cos((2m(N-1)-k)t + q\frac{\pi}{2})\right].$$

The coefficients of the even trigonometric polynomial interpolating the function $f(x)$ on the grid $\Delta 2_N^{(1)}$ are calculated using the formulas:

$$a0 = \frac{2}{N}\sum_{j=1}^{N} f_j, \tag{14}$$

$$a_k = \frac{2}{N}\sum_{j=1}^{N} f_j \cdot \cos(kx_j), \tag{15}$$

and the trigonometric spline $Stc(0,1,\sigma 1,r,q,t)$ and its derivatives of the order $q$ $(q=1,2,...,r)$ are calculated by the formulas

$$Stc(0,1,\sigma 1,r,q,t) = \frac{a_0}{2} I(q) + \sum_{k=1}^{N-1} a_k \frac{C(\sigma 1,r,q,k,t)}{H(\sigma 1,r,k)}, \tag{16}$$

where $I(q)$ is determined by formula (5);

$$H(\sigma 1, r, k) = \sigma 1(r, k) + \sum_{m=1}^{\infty} (-1)^m\left[\sigma 1(r, 2mN+k) + \sigma 1(r, 2mN-k)\right]; \tag{17}$$

$$C(\sigma 1, r, q, k, t) = \sigma 1(r,k)k^q \cos(kt + q\frac{\pi}{2}) +$$
$$+ \sum_{m=1}^{\infty}\left[\sigma 1(r, 2mN+k)(2mN+k)^q \cos((2mN+k)t + q\frac{\pi}{2}) + \right. \tag{18}$$
$$\left. + \sigma 1(r, 2mN-k)(2mN-k)^q \cos((2mN-k)t + q\frac{\pi}{2})\right].$$

Now let's move on to the consideration of odd splines, which we will consider on the grids $\Delta 3_N^{(I)}$ given on $[0, \pi]$.

The coefficients of the odd trigonometric polynomial interpolating the function $f(x)$ at the nodes of the grid $\Delta 3_N^{(0)}$ are calculated using the formulas:

$$b_k = \frac{2}{N+1} \sum_{j=1}^{N} f_j \cdot \sin(kx_j) \quad (19)$$

and the odd spline $Sts(0,0,\sigma1,r,q,t)$ and its derivatives of the order $q$ $(q=1,2,...,r)$, are calculated by the formulas

$$Sts(0,0,\sigma1,r,q,t) = \sum_{k=1}^{N} b_k \frac{S(\sigma1,r,q,k,t)}{H(\sigma1,r,k)}, \quad (20)$$

where

$$H(\sigma1,r,k) = \sigma1(r,k) + \sum_{m=1}^{\infty} \left[\sigma1(r,2m(N+1)+k) + \sigma1(r,2m(N+1)-k)\right], \quad (21)$$

$$S(\sigma1,r,q,k,t) = \sigma1(r,k)k^q \sin(kt + q\frac{\pi}{2}) +$$

$$+ \sum_{m=1}^{\infty} \left[ \sigma1(r,2m(N+1)+k)(2m(N+1)+k)^q \cos((2m(N+1)+k)t + q\frac{\pi}{2}) - \right. \quad (22)$$

$$\left. - \sigma1((r,2m(N+1)-k))(2m(N+1)-k)^q \sin((2m(N+1)-k)t + q\frac{\pi}{2}) \right].$$

The coefficients of the odd trigonometric polynomial interpolating the function on the grid $\Delta 3_N^{(1)}$ are as follows

$$b_k = \frac{2}{N} \sum_{j=1}^{N} f_j \cdot \sin(kx_j) \quad (23)$$

and the trigonometric spline $Sts(1,1,\sigma1,r,q,t)$ and its derivatives of the order $q$ $(q=1,2,...,r)$, are calculated by the formulas

$$Sts(1,1,\sigma1,r,q,t) = \sum_{k=1}^{N-1} b_k \frac{S(\sigma1,r,q,k,t)}{H(\sigma1,r,k)} + .5 b_N \frac{S(\sigma1,r,q,k,t)}{H(\sigma1,r,N)}, \quad (24)$$

$$H(\sigma1,r,k) = \sigma1(r,k) + \sum_{m=1}^{\infty} \left[\sigma1(r,2mN+k) + \sigma1(r,2mN-k)\right] \quad (25)$$

$$S(\sigma1,r,q,k,t) = \sigma1(r,k)k^q \sin(kt + q\frac{\pi}{2}) +$$

$$+ \sum_{m=1}^{\infty} (-1)^m \left[ \sigma1(r,2mN+k)(2mN+k)^q \sin((2mN+k)t + q\frac{\pi}{2}) - \right. \quad (26)$$

$$\left. - \sigma1(r,2mN-k)(2mN-k)^q \sin((2mN-k)t + q\frac{\pi}{2}) \right].$$

Note that the odd splines $Sts(0,0,\sigma,r,t)$ and $Sts(1,1,\sigma,r,t)$ and their even-order derivatives are zero at points $0$ and $\pi$.

Note [18], [19] that the functions $\frac{C(\sigma1,r,q,k,t)}{H(\sigma1,r,k)}$, calculated by formulas (12), (13) and (17), (18), form a system of functions orthogonal to the set of discrete equidistant points on the grids $\Delta 2_N^{(0)}$ and $\Delta 2_N^{(1)}$ respectively, and belong to the class $C_{[0,\pi]}^{r-1}$. With this approach, the coefficients $a_0^{(I)}$, $a_k^{(I)}$ $(k=1,2,...,N)$, can be considered as the coefficients of discrete (finite) pairwise Fourier series. This also applies to the functions $\frac{S(\sigma1,r,q,k,t)}{H(\sigma1,r,k)}$, which are calculated by formulas (21), (22) and (25), (26); these functions form a system of functions orthogonal to the set of discrete equidistant points of the grids $\Delta 3_N^{(0)}$ and $\Delta 3_N^{(1)}$ respectively, and belong to the class $C_{[0,\pi]}^{r-1}$; the coefficients $b_k^{(I)}$ ($k=1,2,...,N$) can also be considered as coefficients of discrete (finite) odd Fourier series.

**3. Representation by trigonometric B-splines.** As is well known (see, e.g., [3] [4]), B-splines play an important role in the theory of polynomial splines, which are the basis in classes of polynomial splines. In [21] [22], the classes of trigonometric B-splines and the classes of kernels of interpolation trigonometric splines were introduced. Trigonometric B-splines are the carriers of information about the smoothness of the spline, and the kernels of trigonometric interpolation splines are the carriers of information about the discrete values of the interpolated function.

By using trigonometric B-splines and kernels of trigonometric interpolation splines, trigonometric interpolation splines can be represented as a convolution of a B-spline of a certain order with a corresponding kernel of a trigonometric interpolation spline. We will not delve into these issues here, referring interested specialists to [18] [19]. We will only give the representation and graphs of trigonometric B-splines for the case of full trigonometric splines discussed in Section 3; we will consider both the Riemann convergence factors $\sigma(r,k)$ and the convergence factors $\sigma 1(r,k)$.

The B-splines with the convergence factor $\sigma 1(r,k)$ of the first and second kind and their derivatives of order $q$ ($q=1,2,...,r$) are calculated by the following formulas

$$BC(r,j,t) = \frac{1}{\pi}\left[\frac{1}{2}I(q) + \sum_{k=1}^{n} C(0,\sigma 1,r,q,k,t-x_j)\right], \qquad (27)$$

$$BC0(r,j,t) = \frac{1}{\pi}\left[\frac{1}{2}I(q) + \sum_{k=1}^{n} \frac{C(0,\sigma 1,r,q,k,t-x_j)}{H(0,0,1+r,k)}\right], \qquad (28)$$

$$BC1(r,j,t) = \frac{1}{\pi}\left[\frac{1}{2}I(q) + \sum_{k=1}^{n} \frac{C(0,\sigma 1,r,q,k,t-x_j)}{H(1,0,1+r,k)}\right], \qquad (29)$$

where the expressions $H(0,0,1+r,k)$ and $C(I1,\sigma 1,r,q,k,t)$ are calculated by formulas (7), (8).

Note that the expressions for calculating the B-splines $BR(r,j,t)$, $BR0(r,j,t)$ and $BR1(r,j,t)$ with the Riemann convergence factor (1) are given in [19].

The graphs of trigonometric B-splines with convergence factors (1) and (2) are shown in Figures 1-3.

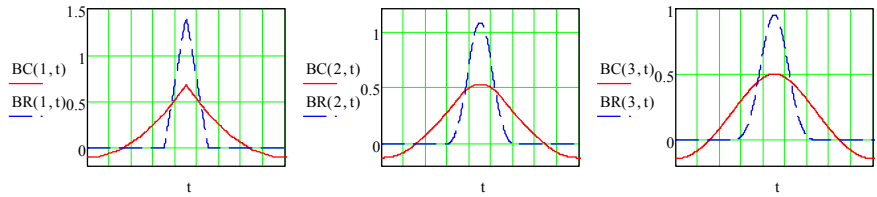

Fig. 1. Trigonometric $BC$- and $BR$-splines of the first kind ($r=1,2,3$)

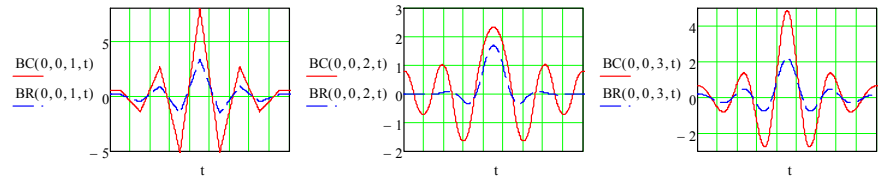

Fig. 2. Trigonometric $BC0$- and $BR0$-splines of the second kind on the same-named grids
$I1=0$, $I2=0$, ($I1=1$, $I2=1$)

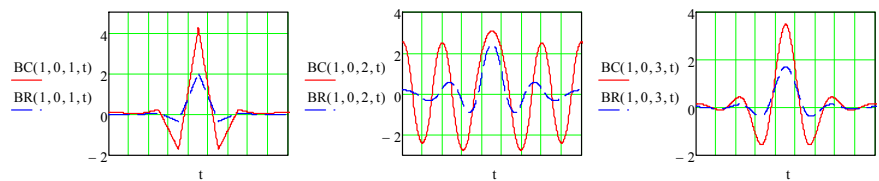

Fig. 3. Trigonometric $BC1$- and $BR1$-splines of the second kind on differently named grids
$I1=1$, $I2=0$, ($I1=0$, $I2=1$)

It is easy to check that the trigonometric B-spline $BC(r, j, t)$ satisfies the following conditions

$$\int_{-\pi}^{\pi} BC(r, j, t) dt = 1 \quad \text{for all } r \text{ and } j. \tag{30}$$

Note that all of the following trigonometric splines $BR(r, j, t)$, $BC0(r, j, t)$, $BR0(r, j, t)$, $BC1(r, j, t)$, $BR1(r, j, t)$ satisfy the same conditions; also, all of these splines belong to the class $C_{[0, 2\pi]}^{r-1}$.

Note that the trigonometric B-splines $BR(r, j, t)$ coincide with polynomial B-splines of the same degree, constructed on uniform grids.

It is known [3], [4] that polynomial interpolation simple splines $S_{2m-1}(r, t)$ can be represented by polynomial B-splines as follows

$$S_{2m-1}(f, \Delta_N, x) = \sum_{j=1}^{N} \alpha_j B_j(2m-1, x), \tag{31}$$

where $\Delta_N$ is an interpolation grid containing $N$ nodes, $f$ is an interpolated function, $B_j(2m-1, x)$ is a polynomial B-spline of degree $2m-1$, $\alpha_j$ are coefficients determined from the interpolation conditions.

A similar representation holds for trigonometric interpolation splines. Indeed, it is known (see, e.g., [17]) that in order for the interpolation problem

$$\sum_{j=1}^{N} \alpha_j \varphi_j(x_k) = f_k, \quad (k = 1, 2, \ldots, N), \tag{32}$$

has a solution and is unique for any real function $f(x)$, it is necessary and sufficient that the determinant

$$D(\varphi_1, \varphi_2, \ldots, \varphi_N) = \begin{vmatrix} \varphi_1(x_1), \varphi_2(x_1), \ldots, \varphi_N(x_1) \\ \varphi_1(x_2), \varphi_2(x_2), \ldots, \varphi_N(x_2) \\ \ldots\ldots\ldots\ldots\ldots\ldots\ldots\ldots\ldots\ldots \\ \varphi_1(x_N), \varphi_2(x_N), \ldots, \varphi_N(x_N) \end{vmatrix} \tag{33}$$

was different from 0 when all $x_1, x_2, \ldots, x_N$ are unequal in pairs.

It is easy to verify that the determinants (33) formed from the trigonometric B-splines $BR(r, j, t)$, $BC(r, j, t)$, $BC0(r, j, t)$, $BR0(r, j, t)$, $BC1(r, j, t)$, $BR1(r, j, t)$ on the grid $\Delta 1_N^{(0)}$, are different from 0 at least for $r = 1, 2, 3, 4, 5, 11$ (except for the determinant $BC(11, j, t)$). Table 1 shows the values of the determinants for $r = 1, 2, 3, 4, 5, 11$ and $N = 9$.

TABLE 1. The values of the determinants (33) for $r = 1, 2, 3, 4, 5, 11$ and $N = 9$

| B-spline type | r=1 | r=2 | r=3 | r=4 | r=5 | r=11 |
| --- | --- | --- | --- | --- | --- | --- |
| $BR(r, j, t)$ | 25.1548 | 1.46797 | 0.3538 | 0.0770 | 0.0189 | 5.893*10^-6 |
| $BC(r, j, t)$ | 9.88*10^-4 | 2.44*10^-8 | 5.5*10^-10 | 1.7*10^-13 | 1.1*10^-15 | 0 |
| $BR0(r, j, t)$ | 6.4396*10^3 | 105.3279 | 512.0283 | 103.5795 | 246.0022 | 117.3284 |
| $BC0(r, j, t)$ | 1.0271*10^6 | 1.1347*10^3 | 8.1665*10^4 | 3.7613*10^3 | 3.9236*10^4 | 1.8713*10^4 |
| $BR1(r, j, t)$ | 434.9783 | 837.8267 | 116.5782 | 324.7136 | 101.8178 | 94.1006 |
| $BC1(r, j, t)$ | 6.9376*10^4 | 9.0262*10^3 | 1.8593*10^4 | 1.1791*10^4 | 1.6239*10^4 | 1.5008*10^4 |

It is worth noting that similar results are observed on the grid $\Delta 1_N^{(1)}$.

From formula (32) it follows that the trigonometric B-splines $BC(r, j, t)$, $BR(r, j, t)$, $BR0(r, j, t)$, $BC0(r, j, t)$, $BC1(r, j, t)$, $BR1(r, j, t)$ form the basis in the space of trigonometric interpolation $2\pi$-periodic splines on the grids $\Delta 1_N^{(l)}$. Thus, the trigonometric interpolation splines on these grids can be represented as

$$St(0, 0, \sigma 1, r, q, x) = \sum_{j=1}^{N} \alpha_j BC0(r, q, x - x_j^{(0)}), \tag{32}$$

where the coefficients $\alpha_j$ are determined from the interpolation conditions. Similar representations are available for the trigonometric B-splines $BC(r,j,t)$, $BR(r,j,t)$, $BR0(r,j,t)$, $BC1(r,j,t)$, $BR1(r,j,t)$. Note that for the B-splines $BC(r,j,t)$ such a representation takes place for small values of the parameter $r$, which is determined by the bit grid of the computer used.

It should be noted that formula (32) has different variants, which are determined by the given stitching and interpolation grids; however, we will not go into this issue here.

**4. Representation through fundamental trigonometric splines.** Let the grid $\Delta_N = \{x_j\}_{j=1}^N$, $x_1 < x_2 < ... < x_N$ be given on the segment $[a,b]$. The system of functions $\varphi_k(x)$, $x \in [a,b]$, is called fundamental on the grid $\Delta_N$ [17] if

$$\varphi_k(x_j) = \begin{cases} 0, & \text{if } k \neq j; \\ 1, & \text{if } k = j. \end{cases} \quad (33)$$

The systems of fundamental functions include, in particular, the system of Lagrange interpolation polynomials, and on uniform grids, the systems of complete trigonometric interpolation polynomials [23]; the same paper also considered systems of fundamental trigonometric splines.

Since it is very convenient to use systems of fundamental functions, we construct systems of fundamental splines on the uniform grids considered above $\Delta 1_N^{(I)}$, $\Delta 2_N^{(I)}$, $\Delta 3_N^{(I)}$.

**Complete fundamental splines on the grids $\Delta 1_N^{(I)}$.** Such splines and their derivatives, which we will denote by $St^*(I1,I2,\sigma1,r,q,t)$ on the grids $\Delta 1_N^{(I)}$ defined on $[0, 2\pi]$, are constructed by the following formulas

$$St^*(I1,I2,\sigma1,r,q,k,t) = \frac{1}{N}\left[I(q) + 2\sum_{j=1}^{n}\left[\frac{C(I1,I2,r,q,k,j,t)}{H(I1,I2,\sigma1,r,j)}\right]\right], \quad (34)$$

where $H(I1,I2,\sigma1,r,k)$ is calculated using formulas (7);

$$C(I1,I2,r,q,k,j,t) = \sigma1(r,j)k^q \cos(j(t-(x_k^{(0)})^{1-I2}(x_k^{(1)})^{I2}) + q\frac{\pi}{2}) +$$

$$+ \sum_{m=1}^{\infty}\left[(-1)^{m\cdot(I1+I2)}\left[\sigma1(r,mN+j)(mN+j)^q \cos\left((mN+j)(t-(x_k^{(0)})^{1-I2}(x_k^{(1)})^{I2}) + q\frac{\pi}{2}\right) + \right. \right. \quad (35)$$

$$\left. \left. +(-1)^{1+r}\sigma1(r,mN-j)(mN-j)^q \cos\left((mN-j)(t-(x_k^{(0)})^{1-I2}(x_k^{(1)})^{I2}) + q\frac{\pi}{2}\right)\right]\right].$$

The fundamental splines of $St^*(I1,I2,\sigma1,r,q,k,t)$ belong to the class $C_{[0,2\pi)}^{r-1}$.

Given the fundamental splines $St^*(I1,I2,\sigma1,r,q,k,t)$, the interpolation trigonometric spline $St(I1,I2,\sigma1,r,q,t)$ and its derivatives of the order $q$ can be represented as follows

$$St(I1,I2,\sigma1,r,q,t) = \sum_{k=1}^{N} f_k^{(I2)} St^*(I1,I2,\sigma1,r,q,k,t), \quad (36)$$

where $f_k^{(I2)}$ ($k=1,2,...,N$) are the value of the interpolated function in the grid nodes $\Delta 1_N^{(I2)}$.

**Incomplete fundamental splines.** Here are the expressions for calculating even $Stc^*(0,0,\sigma1,r,q,k,t)$, $Stc^*(0,1,\sigma1,r,q,k,t)$ and odd $Sts^*(0,0,\sigma1,r,q,k,t)$, $Sts^*(1,1,\sigma1,r,q,k,t)$ ($k=1,2,...,N$) fundamental splines.

The even fundamental splines $Stc^*(0,0,\sigma1,r,q,t)$ and their derivatives on the grid $\Delta 2_N^{(0)}$ are calculated by the expressions

$$Stc^*(0,0,\sigma1,r,q,k,t) =$$

$$= \frac{2}{N-1}\left[.5 \cdot I(q) + \sum_{j=1}^{N-2}\frac{C(\sigma1,r,q,j,t)\cos(x_k^{(0)})}{H(\sigma1,r,j)} + .5\frac{C(\sigma1,r,q,N-1,t)\cos(x_k^{(0)})}{H(\sigma1,r,N-1)}\right], \quad (37)$$

where $H(\sigma 1, r, j)$ is calculated by formula (12);

$$C(\sigma 1, r, q, j, t) = \sigma 1(r, j) j^q \cos(jt + q\frac{\pi}{2}) +$$
$$+ \sum_{m=1}^{\infty} \left[ \sigma 1(r, 2m(N-1) + j)(2m(N-1) + j)^q \cos((2m(N-1) + j)t + q\frac{\pi}{2}) + \right. \tag{38}$$
$$\left. + \sigma 1(r, 2m(N-1) - j)(2m(N-1) - j)^q \cos((2m(N-1) - j)t + q\frac{\pi}{2}) \right].$$

The even fundamental splines $Stc^*(0,1,\sigma 1, r, q, t)$ and their derivatives on the grid $\Delta 2_N^{(1)}$ are calculated by the expressions

$$Stc^*(0,1,\sigma 1, r, q, k, t) = \frac{2}{N} \left[ .5 \cdot I(q) + \sum_{j=1}^{N-1} \frac{C(\sigma 1, r, q, j, t) \cos(j x_k^{(1)})}{H(\sigma 1, r, j)} \right], \tag{39}$$

where $H(\sigma 1, r, j)$ is calculated by formula (17);

$$C(\sigma 1, r, q, j, t) = \sigma 1(r, j) j^q \cos(jt + q\frac{\pi}{2}) +$$
$$+ \sum_{m=1}^{\infty} \left[ \sigma 1((r, 2mN + j)(2mN + j)^q \cos((2mN + j)t + q\frac{\pi}{2}) + \right. \tag{40}$$
$$\left. + \sigma 1(r, 2mN - j)(2mN - j)^q \cos((2mN - j)t + q\frac{\pi}{2}) \right].$$

Given the fundamental splines $Stc^*(0, I2, \sigma 1, r, q, k, t)$, the interpolation trigonometric spline $Stc(0, I2, \sigma 1, r, q, t)$ and its derivatives of the order $q$ can be represented as follows

$$Stc(I1, I2, \sigma 1, r, q, t) = \sum_{k=1}^{N} f_k^{(I2)} Stc^*(I1, I2, \sigma 1, r, q, k, t), \tag{41}$$

where $f_k^{(I2)}$ ($k = 1, 2, ..., N$) are the values of the interpolated function in the grid nodes $\Delta 2_N^{(I2)}$.

The odd fundamental splines $Sts^*(0, 0, \sigma 1, r, q, t)$ on the grid $\Delta 3_N^{(0)}$ are calculated by the expressions

$$Sts^*(0, 0, \sigma 1, r, q, k, t) = \frac{2}{N+1} \sum_{j=1}^{N} \frac{S(\sigma 1, r, q, j, t) \sin(j x_k^{(0)})}{H(\sigma 1, r, j)}, \tag{42}$$

where $H(\sigma 1, r, j)$ is calculated by formula (21);

$$S(\sigma 1, r, q, j, t) = \sigma 1(r, j) j^q \sin(jt + q\frac{\pi}{2}) +$$
$$+ \sum_{m=1}^{\infty} (-1)^m \left[ \sigma 1((r, 2m(N+1) + j))((2m(N+1) + j)^q \sin((2m(N+1) + j)t + q\frac{\pi}{2}) - \right. \tag{43}$$
$$\left. - \sigma 1((r, 2m(N+1) - j))((2m(N+1) - j)^q \sin((2m(N+1) - j)t + q\frac{\pi}{2}) \right].$$

The odd fundamental splines $Sts^*(1, 1, \sigma 1, r, q, k, t)$ on the grid $\Delta 3_N^{(1)}$ are calculated by the expressions

$$Sts^*(1, 1, \sigma 1, r, q, k, t) = \frac{2}{N} \left[ \sum_{j=1}^{N-1} \frac{S1(\sigma 1, r, q, j, t) \sin(j x_k^{(0)})}{H(\sigma 1, r, j)} + .5 \frac{S1(\sigma 1, r, q, N, t) \sin(N x_k^{(0)})}{H(\sigma 1, r, N)} \right], \tag{44}$$

where $H(\sigma 1, r, k)$ is calculated by formula (25);

$$S1(\sigma 1, r, q, j, t) = \sigma 1(r, j) j^q \sin(jt + q\frac{\pi}{2}) +$$
$$+ \sum_{m=1}^{\infty} (-1)^m \left[ \sigma 1((r, 2mN + j))((2mN + j)^q \sin((2mN + j)t + q\frac{\pi}{2}) - \right. \tag{45}$$

$$-\sigma 1((r, 2mN - j))((2mN - j)^q \sin((2mN - j)t + q\frac{\pi}{2})\Big].$$

It is clear that the odd splines $Sts^*(0,0,\sigma 1,r,q,k,t)$, $Sts^*(1,1,\sigma 1,r,q,k,t)$ and their derivatives of odd orders take zero values at the ends of the segment $[0,\pi]$.

Given the fundamental splines $Sts^*(I1,I2,\sigma 1,r,q,k,t)$, the interpolation trigonometric spline $Sts(I1,I2,\sigma 1,r,q,t)$ and its derivatives of the order $q$ can be represented as follows

$$Sts(I1,I2,\sigma 1,r,q,t) = \sum_{k=1}^{N} f_k^{(I2)} Sts^*(I1,I2,\sigma 1,r,q,k,t), \qquad (46)$$

where $f_k^{(I2)}$ ($k = 1,2,...,N$) are the values of the interpolated function in the nodes of the grid $\Delta 3_N^{(I2)}$.

Note that even and odd splines belong to the class $C_{[0,\pi]}^{r-1}$.

Finally, we note that the fundamental splines form orthogonal systems of functions in the sense of a discrete scalar product [17] on the corresponding grids.

## Conclusions

1. A representation of trigonometric $2\pi$-periodic interpolation splines of arbitrary order on $\Delta 1_N^{(I)}$ grids in terms of the coefficients of interpolation trigonometric polynomials on the same grids is given.

2. The functions $\dfrac{C(I1,r,q,k,t)}{H(I1,I2,\sigma 1,r,k)}$ and $\dfrac{S(I1,r,q,k,t)}{H(I1,I2,\sigma 1,r,k)}$, calculated by formulas (7)-(9), belong to the class $C_{[0,2\pi)}^{r-1}$; in addition, these functions form systems of functions orthogonal to the set of discrete equidistant points of the grids $\Delta 1_N^{(I)}$. With this approach, the coefficients of $a_0^{(I)}$, $a_k^{(I)}$, $b_k^{(I)}$ ($k = 1,2,...,n$) can be considered as coefficients of discrete (finite) Fourier series.

3. In the case of interpolation of non-periodic functions, it is advisable to use even or odd trigonometric interpolation splines on the grids $\Delta 2_N^{(I)}$ and $\Delta 3_N^{(I)}$, which are constructed through the coefficients of trigonometric even or odd interpolation trigonometric polynomials on the same grids respectively.

4. The functions $\dfrac{C(\sigma 1,r,q,k,t)}{H(\sigma 1,r,k)}$, calculated by formulas (12), (13) and (17), (18), belong to the class $C_{[0,\pi]}^{r-1}$ and form a system of functions orthogonal to the set of discrete equidistant points of the grids $\Delta 2_N^{(0)}$ and $\Delta 2_N^{(1)}$ respectively. This also applies to the functions $\dfrac{S(\sigma 1,r,q,k,t)}{H(\sigma 1,r,k)}$, calculated by formulas (21), (22) and (25), (26), which belong to the class $C_{[0,\pi]}^{r-1}$ and form a system of functions orthogonal to the set of discrete equidistant points of the grids $\Delta 3_N^{(0)}$ and $\Delta 3_N^{(1)}$ respectively; with this approach, the coefficients $a_0^{(I)}$, $a_k^{(I)}$ ($k = 1,2,...,N-1$) and ($b_k^{(I)}$ $k = 1,2,...,N$) can be considered as coefficients of discrete (finite) even and odd Fourier series.

5. The representation of trigonometric $2\pi$-periodic interpolation splines of arbitrary orders through trigonometric B-splines $BC(r,j,t)$, $BR(r,j,t)$, $BR0(r,j,t)$, $BC0(r,j,t)$, $BC1(r,j,t)$, $BR1(r,j,t)$ of the same orders belonging to the class $C_{[0,2\pi)}^{r-1}$ is investigated. It is shown that these B-splines form a basis in the space of trigonometric $2\pi$-periodic interpolation splines.

6. A representation of trigonometric $2\pi$-periodic interpolation splines of arbitrary orders on $\Delta 1_N^{(I)}$ grids through fundamental trigonometric splines of the same orders on the same grids is given;

such splines belong to the class $C^{r-1}_{[0,2\pi)}$ and form a system of functions orthogonal to a discrete set of points of the grid $\Delta 1^{(I)}_N$.

7. In the case of interpolation of non-periodic functions, it is advisable to use even or odd fundamental trigonometric splines of arbitrary orders on the grids $\Delta 2^{(I)}_N$ and $\Delta 3^{(I)}_N$; such splines belong to the class $C^{r-1}_{[0,\pi]}$ and form a system of functions orthogonal to the discrete set of points of the grids $\Delta 2^{(I)}_N$ and $\Delta 3^{(I)}_N$ respectively.
8. The coefficients in the representations of interpolation splines through fundamental splines are the values of the interpolated function at the nodes of the interpolation grids.
9. The presented representations of trigonometric splines can be used in many problems of both theoretical and applied nature.

**References**

1. Blatter C. Wavelets-Eine Einführung. Friedr. Vieveg & Sohr. Verlagsgesellshaft mbH. Braunschweig/Wiesbaden, 1998.
2. Ahlberg J.H., Nilson E.N., Walsh J.L. The Theory of Splines and Their Applications. – Academic Press. New York/London, 1967.
3. De Boor Carl. A Practical Guide to Splines, Revised Edition. Springer-Verlag, 2003.
4. Завьялов Ю.С., Квасов Б.И., Мирошниченко В.Л. Методы сплайн-функций. – М., Наука, 1980. – 352 с.
5. Корнейчук Н.П. Сплайны в теории приближения. – М.: Наука, 1984. – 356 с.
6. Корнейчук Н.П. Экстремальные задачи теории приближения. – М.: Наука, 1976. – 320 с.
7. Fairweather G., Meade D. A survey of spline collocation methods for the numerical solution of differential equations. – CRC Press, Mathematics for large scale computing, 2020. – 45 p.
8. Chegolin A.P. On the application of the collocation method in spaces of p-summable functions to the Fredgolm integral education of the second kind. – Bulletin of higher educational institutions. North Caucsus Region. Natural Science, 2020. No. 4.
9. Denysiuk V. Generalized Trigonometric Functions and Their Applications // IOSR Journal of Mathematics (IOSR-JM). Volume 14, Issue 6 Ver. I, 2018, pp.19-25.
10. Denysiuk V. About one method of construction of interpolation trigonometric splines. arXiv:1902.07970. –2019.
11. Denysiuk V. Polynomial and trigonometric splines. arXiv:2110.04781.
12. Denysiuk V. Trigonometric splines and some of their properties. arxiv.org/abs/2011.06310.
13. Денисюк В.П. . Тригонометричні ряди та сплайни. – К.:НАУ, 2017. – 212 с.
14. Zygmund A. Trigonometric series. – Volume 1,2, Cambridge at The University Press, 1959.
15. Hardy G.H. Divergent Series. – Oxford, 1949.
16. Бари Н.К. Тригонометрические ряды. – М., Физматгиз, 1961 г. – 936 с.
17. Гутер Р.С., Кудрявцев Л.Д., Левитан Б.М. Элементы теории функций. – М.: Гос. изд-во физ.-мат. лит-ры, 1963 г. – 243 с.
18. Hamming R.W. Numerical Methods for Scientists and Engineers. Second Edition, Dover Publications, Inc., New York, 1973.
19. Hamming R.W. Digital Filtres. – Prentis-Hall, Inc., Englewood Cliffs, New Gersey, 1977.
20. Lanczos C. Applied Analysis. – Prentis-Hall, Inc., 1956.
21. Denysiuk V., Hryshko O. B-splines and kernels of trigonometric interpolation splines arXiv:2310.18722.
22. Denysiuk V., Hryshko O. Generalization of Trigonometric B-splines and Kernels of Interpolating Trigonometric Splines with Riemann Convergence Multipliers arXiv:2310.18722.
23. Denysiuk V. Fundamental trigonometric interpolation and approximating polynomials and splines. arXiv:1912.01918.